\documentclass[10pt,reqno]{article}
\usepackage{amsmath, amsthm, amssymb, stmaryrd}

 \usepackage[T1]{fontenc} 

\DeclareFontFamily{U}{rsfso}{}
\DeclareFontShape{U}{rsfso}{m}{n}{<-> rsfso10}{}
\DeclareMathAlphabet{\mathloop}{U}{rsfso}{m}{n}

\usepackage{hyperref}
\usepackage{enumerate}
\usepackage{url}
\usepackage{color}
\usepackage{tikz}
\usetikzlibrary{patterns}
\usepackage[T1,T2A]{fontenc}
\usepackage[utf8]{inputenc}

\theoremstyle{definition}

\theoremstyle{remark}

\def \le {\leqslant}
\def \ge {\geqslant}

\usepackage{tikz} 
\usetikzlibrary{calc,arrows,decorations.pathreplacing,fadings,3d,positioning}

\usepackage[mathscr]{euscript}

\usepackage[utf8]{inputenc}

\medskip

\begin{document}

 \begin{Huge}
\centerline{On eutactic  forms}
\end{Huge}

\vskip+0.3cm
 \begin{Large}

\centerline{Nikolay Moshchevitin}
\end{Large}

\vskip+1cm

\begin{small}
{\bf Abstract:}
 We prove two  new  
 results about the number of classes of positive definite eutactic forms in $n$ variables.
  
 \end{small}

\vskip+1cm
{\bf 1. Extremality, eutaxy and perfectness.}
\vskip+0.3cm

In this paper we deal with positive definite quadratic forms
\begin{equation}\label{ququ}
Q({\bf x}) = {\bf x}^\top \frak{Q}{\bf x} =
\sum_{i,j=1}^n q_{i,j} x_ix_j,\,\,\,\,\, q_{i,j} = q_{j,i}
\end{equation}
  in $n$ real variables $ x_1,...,x_n$. Here $\frak{Q}$ is a symmetric $n\times n$ matrix. We consider  $\frak{Q}$ as a point in Euclidean space $\mathbb{R}^N, N = \frac{n(n+1)}{2}$
  with  coordinates $q_{i,j} 1\le i\le j\le n$ and inner product
  $$
  \langle \frak{Q}, \frak{Q}'\rangle = {\rm tr}\, \frak{Q}\frak{Q}' = 
  \sum_{i,j=1}^n q_{i,j} q_{i,j}'.
  $$ 
Recall that two  forms
$Q'({\bf x})$ and $Q''({\bf x})$
are called {\it (integer) equivalent } if there exists a   $n\times  n$ matrix $\frak{U}$ with  integer entries and determinant
${\rm det}\, \frak{U} = \pm 1$ such that 
$Q''({\bf x})=Q'(\frak{U}{\bf x})$ for all ${\bf x}\in \mathbb{Z}^n$.
For any positive definite quadratic form $Q$ in the case $n\ge 11$ for the set of its automorphisms
$$
{\rm Aut} (Q) =
\{ \frak{U} \in {\rm GL}_n (\mathbb{Z}): \,\,\,  \frak{U}^\top  \frak{Q}  \frak{U} =  \frak{Q}\}
$$
we have  the upper bound
\begin{equation}\label{aut}
|{\rm Aut} (Q)|  \le 2^n n!
\end{equation}\
due to Feit–Friedland's theorem \cite{FF}.

 For every positive definite quadratic form $Q$ there exists real $n\times n $ matrix $\frak{A}_Q$   such that
$ \frak{Q} = \frak{A}_Q^\top  \frak{A}_Q  $ so one can associate with $Q$ a lattice $\Lambda_Q = \frak{A}_Q \mathbb{Z}^n$.  
  Consider  the 
arithmetical minimum
 $$
m_Q = \min_{{\bf x} \in \mathbb{Z}^n\setminus\{{\bf 0}\}} Q({\bf x})
$$
of the form $Q({\bf x})$ and the set 
$$
\Sigma_Q = \Sigma_Q^0/{\pm 1},\,\,\,\,\, \Sigma_Q^0 = \{{\bf x}\in \mathbb{Z}^n: \,\,\, Q({\bf x}) = m_Q\} 
$$
of all integer representations of the arithmetical minimum (here we should note that if ${\bf v}$ is a minimal vector, then $-{\bf v}$ is also a minimal vector, so in $
\Sigma_Q$  we consider just one representative for  a pair ${\bf v} , -{\bf v} $).
In other words,  the set
\begin{equation}\label{seet}
S_Q  =S_Q^0 / \{\pm 1\},\,\,\,\,\,
S_Q^0 = \{ {\bf z}\in \mathbb{R}^n: \exists \, {\bf x} \in \mathbb{Z}^n\,\,\text{such that}\,\,
Q({\bf x}) = m_Q\,\, \text{and}\,\,  {\bf z} =\frak{A}  {\bf x} \} 
\end{equation}
is the set of all non-zero vectors of minimal Euclidean length in the lattice 
$\Lambda_Q$, up to sign $\pm$.
For each $ {\bf v} \in \Sigma_Q$ consider $n\times n$ matrix $\frak{W}_{\bf v} = {\bf v}\cdot {\bf v}^\top$.

Also we consider ellipsoid
$$
\mathcal{E}_Q 
 = \bigl\{{\bf x}\in\mathbb{R}^{n}:Q({\bf x})\le m_Q \bigr\}.
 $$
This ellipsoid
does not contain nonzero integer points inside. Integer points on the boundary of 
$\mathcal{E}_{Q}$ are just points $\pm {\bf v},  {\bf v}\in \Sigma_{Q}$.

\vskip+0.3cm
We recall the definitions of {\it extremal form},  {\it eutactic form},
{\it weakly eutactic form}
 and  {\it perfect form}.
 
 \vskip+0.3cm
 
 Form $Q({\bf x})$ is called {\it extremal} if  any determinant-preserving   small perturbation of its coefficients $q_{i,j}$ 
makes the arithmetical minimum of the form  $Q({\bf x})$ smaller. 
%The set of all classes of equivalences of positive definite extremal forms in $n$ variables we denote by $\mathloop{E}_n$.

  \vskip+0.3cm
  Form $Q({\bf x})$ is called {\it eutactic} if   there exists a collection of positive numbers 
  $$
  \lambda_{\bf v}>0,\,\,\,\,\,  {\bf v} \in \Sigma_Q
  $$
  such that for the inverse matrix $\frak{Q}^{-1}$  we have equality
  \begin{equation}\label{euz}
  \frak{Q}^{-1} = \sum_{ {\bf v} \in \Sigma_Q} \,   \lambda_{\bf v} \frak{W}_{\bf v}.
  \end{equation}

  \vskip+0.3cm
  Form $Q({\bf x})$ is called {\it  weakly eutactic} if   there exist  real (not necessarily positive) numbers    $
  \lambda_{\bf v},  {\bf v} \in \Sigma_Q
  $ such that (\ref{euz}) holds, that  is matrix $  \frak{Q}^{-1} $ belongs to the subspace in $\mathbb{R}^N$ spanned by the set 
 $ \{\frak{W}_{\bf v}: {\bf v}\in \Sigma_Q\}$. Obviously, eutactic  form $Q$ is automatically  weakly eutactic. Here we also should notice that 
\begin{equation}\label{wellro}
\text{every  weakly eutactic form}\,\, \, Q\, \,\,\text{  is {\it well rounded}},
\end{equation}
 that is the set $ \Sigma_Q^0$ contains $n$ linearly independent vectors, 
 see for example
 Chapter 3, §3.2 from \cite{Martinus}.

  \vskip+0.3cm

Following Voronoï \cite{V} form $Q({\bf x}) $ is defined to be {\it perfect }
if coefficients $q_{i,j} , 1\le i,j \le n$ are uniquely determined by the value of arithmetic minimum 
$
m_Q  
$
and the set
$
\Sigma_Q 
$,
 that is  given $m$ and $\Sigma$ the form 
$Q({\bf x})$
with the given value of $ m_Q =m$ and with given set $\Sigma_Q =\Sigma$ is unique. In particular this means that the whole set 
of symmetric matrices $\mathbb{R}^N$ is spanned by the set 
 $\{ \frak{W}_{\bf v}:{\bf v}\in \Sigma_Q\}$. It is clear from the definition that every perfect form is weakly eutactic.
 
 \vskip+0.3cm
It is known due to Voronoï \cite{V}  that positive definite quadratic form is extremal if and only if it is simultaneously perfect and eutactic. Also  for any dimension $n$ the set of classes of 
equivalences of perfect forms (Voronoï \cite{V}), eutactic forms (Ash \cite{Ash}) and weakly eutactic forms (Berg\'e, Martinet \cite{BM},  see also Ryshkov \cite{RR} however in \cite{RR} the notion of weak eutaxy is not formulated explicitly) is finite.

\vskip+0.3cm
{\bf 2. Results.}

\vskip+0.3cm

Let $ E_n, EU_n, W_n $ and $P_n$ be the numbers of all non-equivalent  extremal, eutactic, weakly eutatic  and perfect positive quadratic forms in  forms in $n$ variables respectively. 
From the discussion  of the previous section is clear that 
$$
E_n \le \min ( P_n, EU_n) \le \max( P_n, EU_n)\le W_n.
$$
It was recently proven that 
\begin{equation}\label{rese}
\exp\left(\frac{n^2\log n}{2} - O(n^2)\right)\le
E_n \le \exp\left(\frac{5n^2\log n}{2} +O(n^2)\right)
\end{equation}
The lower bound here is due to the author \cite{MOS}, the upper bound is due to van Woerden \cite{W}\footnote{In \cite{W} the author proved an explicit and effective upper bound (Theorem 3.1), the constant $\frac{5}{2}$ in the exponent in the right hand side of (\ref{rese}) was not explicitly stated but it immediately follows from calculations. AI tools  (Claude, Fable 5.1) claim that the constant  $\frac{5}{2}$  may be reduced to  $\frac{9}{4}$  }.

\vskip+0.3cm
Now we formulate the main results of the present paper.
\vskip+0.3cm

{\bf Theorem 1.} {\it
For the  number  $W_n$of  pairwise non-equivalent  weakly eutactic  forms in $n$ variables
we have upper bound
$$
W_n \le
 \exp\left(\frac{5n^2\log n}{2} + O(n^2)\right)
.
$$
 }

\vskip+0.3cm

The proof of Theorem 1 is some kind of modification of van Woerden's proof
\cite{W}  of  the upper bound from (\ref{rese}).
 
 \vskip+0.3cm
 
Our next result provides a lower bound for the number of forms which are perfect but not eutactic.
\vskip+0.3cm

{\bf Theorem 2.} {\it
The number of pairwise non-equivalent perfect forms in $n$ variables which are not
eutactic is
$$ 
\ge \exp\left(\frac{n^2\log n}{2} - O(n^2)\right).
$$
}

\vskip+0.3cm

To prove Theorem 2  we follow the steps of the proof of the main result from  \cite{MOS}. The argument is quite similar so we do not give the very details but explain the construction. The difference with \cite{MOS} is that instead of the family of  first Voronoï perfect forms we use Barnes' family of perfect forms from \cite{BA}.

\vskip+0.3cm
The structure of the rest of the paper is as follows. In Section 3 we recall all the basic facts concerning  geometry  of Ryshkov's and Voronoï's perfect polyhedra.
In Section 4 we explain how we adapt Woerden's  proof to establish the result of Theorem 1. In Section 5 we recall all necessary facts concerning Barnes' family of non-eutactic perfect forms. Section 6 deals with general versions of auxiliary statements from \cite{MOS}. In Section 7 we finalise the proof of Theorem 2; this final argument is close to the proof of Minkowski-Hlawka theorem (see Chapter VI from  \cite{Cas} and similar construction in \cite{corobo}.

\vskip+0.3cm
{\bf 3. Geometry of perfect Ryshkov's polyhedra.}
\vskip+0.3cm

In this section we recall the definition and all necessary property of Ryshkov's polyhedron $\mu (m)$. All these properties are explained and proved in papers 
\cite{R,RR}. By
$$
\mathloop{K} = \{ \frak{Q}\in \mathbb{R}^N:\,\,\, 
Q\,\,\text{is positive definite}\} \subset
\overline{\mathloop{K}} = \{ \frak{Q}\in \mathbb{R}^N:\,\,\, 
Q\,\,\text{is positive semidefinite }\} \subset
 \mathbb{R}^N
$$
we denote the cone of positivity in the space of parameters $\mathbb{R}^N$.

\vskip+0.3cm
Let $ m >0$.
Consider primitive integer point ${\bf x} =(x_1,...,x_n)\in \mathbb{Z}^n$ and the corresponding 
Voronoï  point  $\frak{W}_{\bf x} \in \mathbb{R}^N$.
 We associate to $ {\bf x}$  the half-space $ H_{\bf x}$  of the space $\mathbb{R}^N$ defined as
 $$
  H_{\bf x} = \left\{ \frak{Q}\in \mathbb{R}^N:\,\,\,
 \sum_{i,j}  x_ix_jq_{i,j}\ge m\right\}.
 $$ 
 Consider   the intersection 
 $$
 {\bf M}(m)  = \bigcap_{{\bf x} -\,\,\text{primitive}}  H_{\bf x} \,\,\,\,\,
 \text{and its boundary}\,\,\,\,\,
 \mu (m) = \partial \, {\bf M} (m) .
 $$
 
 We  formulate several basic properties of these sets below.

  \vskip+0.3cm 
\noindent
${\bf 1}^0$. The set  $\mu (m)$ is the locus of the points 
$\frak{Q} \in \mathbb{R}^N$
corresponding to positive quadratic forms
$Q$
 with arithmetic minimum equal to $m$.

\noindent
${\bf 2}^0$. The sets  $\mu (m), {\bf M} (m)$ are invariant under the transformations of the group 
$
\mathloop{G} = \{ \frak{G}\}
$
of linear transfromations
$$
\frak{G}:\,\,\, \mathbb{R}^N\to \mathbb{R}^N,\,\,\,\,\,\,
\frak{Q} \mapsto \frak{G}^\top \frak{Q} \frak{G},\,\,\,\,\, \frak{G} \in {\rm GL}_n(\mathbb{Z}).
$$

\noindent
 ${\bf 3}^0$. The set  ${\bf M}  (m)$ is a locally finite polyhedron. All $(N-1)$-dimensional faces of the polyhedron ${\bf M}  (m)$ are unbounded and are pairwise equivalent under the group $\mathloop{G}$.

\noindent
${\bf 4}^0$. The vertices of the polyhedron $\mu (m)$ are in one-to-one correspondence with perfect forms with arithmetic minimum equal to $m$.  

\noindent
${\bf 5}^0$. Polyhedron $\mu (m)$  has finite faces, that is faces which are $l$-dimensional polyhedra  ($0\le l \le N-n = \frac{n(n-1)}{2}$)  with finite number of vertices.
The bound $l\le N-n$  follows from \cite{RR}, \S1 ($2^0$), Definition  and Theorem 2  (every finite face is a face of a complete finite face, which has dimension $N-n$ in terminology of \cite{RR}).

\vskip+0.3cm
The next property  which is in fact Theorem 1 from \cite{RR} is especially important for us.

\vskip+0.3cm
\noindent
${\bf 6}^0$. Quadratic form $Q$ with minimum $m = m_Q$ is well rounded if and only if there exists a {\it finite} face $\mathloop{F}$  of
polyhedron $\mu (m) $ such that $ \frak{Q}\in   \mathloop{F}$.
 \vskip+0.3cm

We now indicate  connection between the polyhedron $\mu (m)$ and Voronoï polyhedron $\pi$ (the perfect partition). First of all we recall the basic definitions and facts.

For a perfect form $Q$ consider cones
\begin{equation}\label{vory}
\pi_Q  =\left\{ \frak{Q}'\in \mathbb{R}^N:\,\,\, \frak{Q}' = \sum_{{\bf v}\in \Sigma_Q} \lambda_{\bf v} \frak{W}_{\bf v} ,\,\,\,\,\,\,  \lambda_{\bf v}> 0 \,\, \forall \, {\bf v} \in \Sigma_Q\right\} \subset\mathloop{K}
\end{equation}
and their closures
$$
\overline{\pi_Q} =\left\{ \frak{Q}'\in \mathbb{R}^N:\,\,\, \frak{Q}' = \sum_{{\bf v}\in \Sigma_Q} \lambda_{\bf v} \frak{W}_{\bf v} ,\,\,\,\,\,\,  \lambda_{\bf v}\ge 0 \,\, \forall \, {\bf v} \in \Sigma_Q\right\}.
$$
The union 
\begin{equation}\label{partai}
\bigcup_{Q -\,\, \text{perfect}}{\pi_Q} \subset\mathloop{K}
\end{equation}
which is a proper subset of $\mathloop{K}$,
corresponds to the partition of the positivity cone  $\mathloop{K}$
known as Voronoï perfect partition.
Voronoï polyhedron 
$$
\pi ={\rm conv}\, \left( \bigcup_{{\bf x}- \,\,\text{primitive integer point}} \frak{W}_{\bf x}\right)
$$ 
satisfies 
$$
\bigcup_{Q -\,\, \text{perfect}}\overline{\pi_Q} =\{ t\frak{X}:\,\,\,  \frak{X}\in \pi, t\ge0\}.
$$
By Voronoï's theorem \cite{V}, 
there exist  non-equivalent perfect forms $Q_1,...,Q_{P_n}$ such that
turnications $ \pi_{Q_i}\cap \mathloop{T}$ are compact and 
for every perfect $Q$ its 
 cone $\pi_Q$ is $\mathloop{G}$-equivalent to one of $ \pi_{Q_i}$.

The next property is 
 Theorem 1 from \cite{R}.

  \vskip+0.3cm
\noindent
${\bf 7}^0$.  Polyhedron $\mu(m)$ is dual to the classical  Voronoï  perfect polyhedron $\pi$, in particular
the number of all $j$-dimensional faces of $\mu(m)$ passing through a fixed vertex $\frak{Q}$ of   $\mu (m) $ is equal to  the number of $N-j$ dimensional faces of the  closed cone 
$\overline{\pi_Q}$.

 \vskip+0.3cm
{\bf Lemma 1.}
{\it The value $W_n$  is bounded from above by the total number of pairwise non-equivalent
finite  faces of polyhedron $\mu (m)$.}
 
  \vskip+0.3cm
  
  Proof.  By (\ref{wellro}) and property ${\bf 6}^0$ we see that any weakly eutactic form $Q$  with mininum $m = m_Q$ corresponds to an inner point 
  $\frak{Q}$ 
  of a unique finite face $\mathloop{F}$ of polyhedron $\mu (m)$  (this face $\mathloop{F}$  contains point $\frak{Q}$ in its relative interior;
 in the case when face is 0-dimensional, that is a vertex of $\mu (m)$, the 
  standard convention is that the relative interior of a point is the point itself).

 Consider  connected equidiscriminant surface  
 $$ \mathloop{D} = \{ \frak{X}\in \mathloop{K}:\, {\rm det}\,( \frak{X})=
   {\rm det}\,( \frak{Q})\} \subset\mathloop{K}$$
     passing through $\frak{Q}$ and the tangent hyperplane $\mathloop{T}$ to $\mathloop{D}$ at point $\frak{Q}$. It follows from weak eutaxy that $ \mathloop{F}\subset \mathloop{T}$, so $\mathloop{F}$ is tangent to  $\mathloop{D}$.
   But it is well known 
   (Minkowski's  determinant inequality,
   see for example \S 13.5 from \cite{Lmi}) that the 
equidiscriminant surface  is strictly convex. 
So in the relative interior of  face $\mathloop{F}$ there may be not more than just one point where $\mathloop{F}$ is tangent to  $\mathloop{D}$.$\Box$

\vskip+0.3cm
{\bf 4. Van Woerden's argument and proof of Theorem 1.}
\vskip+0.3cm 
Let 
$$
\mathloop{T} =
\{ \frak{X}\in \overline{\mathloop{K}}: \,\,   \langle \frak{Q}, \frak{I}\rangle = {\rm tr}\, \frak{Q} \le 1\},\,\,\,\,
\text{(here $\frak{I}$ is $n\times n$ unit matrix})
$$
and   for Jordan-measurable $ \Omega \subset \mathloop{K}$ consider "reduced volume"
$$
V(\Omega) = {\rm Vol}_N  (\Omega \cap \mathloop{T} )
$$
(here 
$\mathrm{Vol}_N$ denotes Lebesgue measure in the coordinates $(q_{i,j})_{i\le j}$).

  We use van Woerden parameters $\ell_n, u_n$ which were introduced in \cite{W}.
  However we do not use here the map $\phi(\cdot)$ so  we recalculate the corresponding volumes and use parameter $\ell_n' = 2^{(n-N)/2}\ell_n.$ The formulas defining Woerden's parameters are as follows:
$$
\ell_n  = 2^{(N-n)/2} \ell_n' = \frac{1}{N!}\frac{2^{(N-n)/2}}{(\frac{1}{8} n^3(n+7))N},\,\,\,\,\,\,\,\,
u_n    = 2^{(N-n)/2} u_n' = \frac{1}{N} \frac{1}{\sqrt{n}} \left(\frac{n-1}{n}\right)^{(N-1)/2}\frac{\pi^{(N-1)/2}}{\Gamma (N/2+1/2)},
$$
\begin{equation}\label{ppa}
 \frac{u_n'}{\ell_n'}  =
 \frac{u_n}{\ell_n} = \exp
\left(\frac{5n^2\log n}{2} + O(n^2)\right)
\end{equation}
 .

The main corollary from  van Woerden's argument \cite{W} is as follows.
Recall that $P_n$ is the number of classes of non-equivalent  positive definite perfect quadratic forms in $n$ variables.
There exist a collection $Q_i, 1\le i \le P_n$  of pairwise non equivalent positive definite perfect quadratic forms
which form a fundamental domain of the perfect partition (\ref{partai}) and 
 such that  in addition
\begin{equation}\label{vory1}
P_n \ell_n' \le  
\sum_{i=1}^{P_n} V(\pi_{Q_i}) \le  V (\overline{\mathloop{K}}) \le u_n'.
\end{equation}
The forms $Q_i$ are chosen to be reduced in the  sense Lemma 3.3 from \cite{W}, in fact this ensures the lower bound 
(\ref{lowerb})  few lines below.

In our consideration each $\pi_{Q_i}$ is partitioned into several simplicial cones 
$\mathloop{S}_{i,\nu} , 1\le \nu \le h_i$ with vertices at certain points of the  form 
$ \frak{W}_{\bf v}$ with ${\bf v} \in \Sigma_{Q_i}$.
The existence of this partition follows from the possibility of triangulation  
of the facets
 $$
 \mathrm{conv}\{\mathfrak W_{\mathbf v}:\mathbf v\in\Sigma_{Q_i}\}
 $$
  of Voronoï's polyhedron $\pi$  (a cross‑section of $\pi_{Q_i}$),
  which is in fact a finite $(N-1)$-dimensional polyhedron,
into $(N-1)$-dimensional simplices with vertices in $\bigcup_{{\bf v}\in \Sigma_{Q_i}}  \frak{W}_{\bf v}$.

For each of these simplicial cones we have  lower bound
\begin{equation}\label{lowerb}
V(\mathloop{S}_{i,\nu} ) \ge 
\ell_n'
\end{equation}
which comes from the argument quite similar to  Corollary 3.4 from \cite{W}.
So instead of van Woerden's inequality (\ref{vory1}) we get 
\begin{equation}\label{vory2}
 \ell_n'\sum_{i=1}^{P_n}  h_i \le  
\sum_{i=1}^{P_n} 
\sum_{\nu = 1}^{h_i} 
V(\mathloop{S}_{i,\nu} )  =
\sum_{i=1}^{P_n} V (\pi_{Q_i})\le 
 V(\overline{\mathloop{K}})
\le u_n'
\end{equation}
(we use the fact that cones $\pi_{Q_i}$ have pairwise disjoint interior).
But  the number of all faces of $\pi_{Q_i}$ is less than the sum of the number of all faces of simplicial cones $\mathloop{S}_{i,\nu} , 1\le \nu \le h_i$.
Each simplicial cone  has $2^N-1$ non-empty faces.
So
the total number of faces of $\pi_{Q_i}$ is not greater than
$$
(2^N-1) h_i
.
$$ 
Every finite face of $\mu (m)$  has a vertex, i.e. a perfect form, which is $\mathloop{G}$‑equivalent to some $Q_i$. 
Hence every finite face is equivalent to a face of  $\mu (m)$  passing through some $Q_i$; by 
${\bf 7}^0$
the faces of  $\mu (m)$ through $Q_i$ correspond bijectively to the non‑zero faces of the cone $\pi_{Q_i}$, whose number is $\le (2^N-1)h_i$.
So
 using Lemma 1 we conclude that 
 $$
 W_n \le (2^N-1) \sum_{i=1}^{P_n} h_i \le  (2^N-1) \frac{u_n'}{\ell_n'}.
$$
Taking into account (\ref{ppa}) we finalise the proof of Theorem 1.$\Box$

  \vskip+0.3cm
{\bf 5. Barnes'  family of forms.}
\vskip+0.3cm

In this section we  formulate  all the necessary results from Barnes'  paper  \cite{BA}.
Consider  the family of positive definite quadratic forms
$$f_n({\bf x})=
f_n(x_1,...,x_n) = \left(\sum_{i=1}^n x_i\right)^2 + \sum_{j=1}^3 \left(x_j^2-x_jx_{j+3}+x_{j+3}^2\right) +\sum_{k=7}^nx_k^2 = {\bf x}^\top \frak{F}_n{\bf x}  ,
$$
where $\frak{F}_{n} = \frak{A}_{f_n}^\top  \frak{A}_{f_n} $ is a  certain $n\times n$ symmetric matrix and
$$
{\rm det}\, \frak{A}_{f_n}  =  \frac{3\sqrt{3(n+7)}}{8},
$$
see  formula (15) from \cite{BA}.
For lattice $ \Lambda_{f_n}$ and its dual 
 $$
  \Lambda_{f_n}^* 
  =
  \{ {\bf x} \in \mathbb{R}^n: \,\,\, ({\bf x},{\bf y}) \in \mathbb{Z}\,\,\,\, \forall\, {\bf y} \in \Lambda_{f_n}\}
  $$
  from formulas (15,16)  from \cite{BA}  and standard description of dual lattice
  we have
  \begin{equation}\label{dudud}
  2\Lambda_{f_n} \subset \Lambda_{f_n},\,\,\,\,\, 3(n+7) \Lambda_{f_n}^* \subset \Lambda_{f_n}.
  \end{equation}

It is known  that for every $ n\ge 11$ the form   $f_n({\bf x})$ is perfect but not extremal and so not eutactic (see \cite{BA}, last statement on p. 154). Moreover,  $m_{f_n} = 2$ and the set $\Sigma_{f_n}$ consists of just $\frac{n(n+1)}{2}$ vectors and so the 
eutaxy coefficients $ \lambda_{\bf v} [f_n] ,{\bf v}\in \Sigma_{f_n}$
are defined from the equation
\begin{equation}\label{EU}
  \frak{F}_n^{-1} = \sum_{ {\bf v} \in \Sigma_{f_n}} \,   \lambda_{\bf v} [f_n] \frak{W}_{\bf v}
\end{equation}
uniquely. In particular, for $n \ge 12$ one of these coefficients is negative.
 For the volume  $\omega_n$ of $\mathcal{E}_{f_n}$
we have equality 
\begin{equation}\label{wo}
\omega_n = 
\operatorname{vol} \mathcal{E}_{f_n}
  = \frac{2^{n/2}\,\pi^{n/2}}{\Gamma\!\left(\tfrac{n}{2}+1\right)|{\det \frak{A}_n}|}
  = \frac{8}{3\sqrt{3}}\cdot
    \frac{(2\pi)^{n/2}}{\sqrt{n+7}\;\Gamma\!\left(\tfrac{n}{2}+1\right)} ,\,\,\,\
\log \omega_n
  = -\frac{n}{2}\log n+ O(n),\,\,\, n \to \infty.
\end{equation}

\vskip+0.3cm
{\bf 6. Statements about lattices.}
\vskip+0.3cm

\vskip+0.3cm
Let  $p>n+1$ be prime and 
 \begin{equation}\label{aa}
 {\bf a}= (a_1,...,a_{n-1}, 1)^\top\in \mathbb{Z}^{n},\,\,\,\,\, 0\le a_j \le p-1, \,\,\ j=1,...,n-1
 \end{equation}
 be an integer vector.
 Define $\mathscr{W}$ to be the set of all $p^{n-1}$ integer vectors of the form (\ref{aa})
 and consider  
\begin{equation}\label{ppil}
 n\times n \,\,\,\,\,\text{matrices}\,\,\,\,\, 
   \frak{B}_{\bf a} = 
 \left(
 \begin{array}{ccccc}
 1&0& \cdots &0&\frac{a_1}{p}\cr
  0&1& \cdots &0&\frac{a_2}{p}\cr
   \vdots&\vdots& \vdots &\vdots&\vdots\cr
    0&0& \cdots &1&\frac{a_{n-1}}{p}\cr
     0&0& \cdots &0&\frac{1}{p}
 \end{array}
 \right),\,\,\,\,\,
 \frak{A}_{{\bf a}}=
 \frak{A} \frak{B}_{\bf a}
 \,\,\,\,\,\,\,\,\,\,
  \text{and lattice} 
  \,\,\,\,\,
    \mathcal {L}_{\bf a} = \frak{A}_{\bf a} \mathbb{Z}^n.
    \end{equation}
Given arbitrary positive definite quadratic form (\ref{ququ})  and integer vector ${\bf a}$ satisfying (\ref{aa}) we consider quadratic form
$$
Q_{\bf a} ({\bf x}) = {\bf x}^\top \frak{Q}_{\bf a}{\bf x}  =
Q( \frak{B}_{\bf a}{\bf x}),\,\,\,\,\,\text{where}\,\,\,\,\,
\frak{Q}_{\bf a}=   \frak{A}_{{\bf a}}^\top  \frak{A}_{{\bf a}} = \frak{B}_{\bf a}^\top\frak{Q}\frak{B}_{\bf a}.
$$

  \vskip+0.3cm
  {\bf Lemma 2.} {\it
  Assume that 
  $$\frak{Q}_{\bf a}=   \frak{A}_{{\bf a}}^\top  \frak{A}_{{\bf a}}\,\,\,\,\,
  \text{with}\,\,\,\,\,
  \frak{A}_{{\bf a}} = \frak{A}_{Q_{\bf a}} =  \frak{A}_{Q} \frak{B}_{{\bf a}}.
  $$
  Assume that 
  for the sets defined in (\ref{seet})  with respect to  lattices 
  $ \frak{A}_{Q} \mathbb{Z}^n$ and   $ \frak{A}_{Q_{{\bf a}}} \mathbb{Z}^n$
  one has 
  \begin{equation}\label{0}
  S_Q = S_{Q_{\bf a}}.
  \end{equation}
   Then
  
    \noindent
  {\rm ({\bf i})} $Q$ is eutactic if and only if  $Q_{\bf a}$ is eutactic;

    \noindent
  {\rm ({\bf ii})} $Q$ is perfect if and only if  $Q_{\bf a}$ is perfect;

  \noindent
  {\rm ({\bf iii})} $Q$ is extremal if and only if  $Q_{\bf a}$ is extremal.

   }

 \vskip+0.3cm
 
 Proof.
     \noindent
  {\rm ({\bf i})}  It follows from the condition (\ref{0}) that 
  $\Sigma_{Q_{\bf a}} = 
  \frak{B}_{\bf a}^{-1} \Sigma_Q$. As 
  $$\frak{Q}_{\bf a}^{-1} =   \frak{B}_{\bf a}^{-1}\frak{Q}_{\bf a}^{-1}  (\frak{B}_{\bf a}^{-1})^\top 
  \,\,\,\,\,
  \text{and}
  \,\,\,\,\,
     \frak{B}_{\bf a}^{-1} \frak{W}_{\bf v}  (\frak{B}_{\bf a}^{-1})^\top  =
     (\frak{B}_{\bf a}^{-1} {\bf v})\, \cdot\, (\frak{B}_{\bf a}^{-1}{\bf v})^\top,
  $$ 
   we see that  equality  (\ref{euz})  is equivalent to
  $$
  \frak{Q}_{\bf a}^{-1}  =\sum_{{\bf v}\in \Sigma_Q} \lambda_{\bf v} \frak{B}_{\bf a}^{-1} \frak{W}_{\bf v}(\frak{B}_{\bf a}^{-1})^\top=
  \sum_{{\bf v}\in \Sigma_Q} \lambda_{\bf v}      (\frak{B}_{\bf a}^{-1} {\bf v})\cdot\, (\frak{B}_{\bf a}^{-1}{\bf v})^\top=
  \sum_{{\bf v}\in \Sigma_{Q_{\bf a}}} \lambda_{\bf v}  \frak{W}_{\bf v}
  $$
  with the same coefficients $\lambda_{\bf v}$.

       \noindent
       Statements
  {\rm ({\bf ii})}  and     
  {\rm ({\bf iii})}  are obvious generalisations of Lemma 3 and Remark after it from \cite{MOS}.$\Box$
      \vskip+0.3cm
          \vskip+0.3cm
          
        We should notice that proof of Theorem 2 uses only   {\rm ({\bf i})}   and   {\rm ({\bf ii})}.     
  
     \vskip+0.3cm
 {\bf Lemma 3.}  
  {\it
  Let $\frak{A}$ be non-degenerate real $n\times n$ matrix and $\frak{G} = \frak{A}^\top\frak{A}$.
  Consider quadratic form $\varphi ({\bf x}) = {\bf x}^\top \frak{G}{\bf x}$ and lattice
  $\Lambda = \frak{A}\mathbb{Z}^n$.
  Let 
  $$
  \Lambda^* 
  =
  \{ {\bf x} \in \mathbb{R}^n: \,\,\, ({\bf x},{\bf y}) \in \mathbb{Z}\,\,\,\, \forall\, {\bf y} \in \Lambda\}
  $$
be the lattice dual to $\Lambda$
  (then in particular  for any integer vector ${\bf x}\in \mathbb{Z}^n$ one has $ M_1\varphi ({\bf x}) \in \mathbb{Z}$).
  
  Assume that 
    \begin{equation}\label{mppi}
  M_1\Lambda \subset \Lambda^*\,\,\,\,\,\text{and}\,\,\,\,\,
  M_2\Lambda^* \subset \Lambda
  \end{equation}
  Consider  prime $p$ such that 
  \begin{equation}\label{ppi}
  M_1M_2 \not\equiv 0 \pmod{p}.
    \end{equation} 
Assume that ${\bf a}\in \mathscr{W}$  satisfies 
          \begin{equation}\label{ppiV}
  M_1\varphi ({\bf a}) \not\equiv 0 \pmod{p},
    \end{equation}

 Assume that a lattice $A$ satisfies
  \begin{equation}\label{la}
 M_2\mathcal{L}_{\bf a}^* \subset A\subset \mathcal{L}_{\bf a} \,\,\,\,\,
 \text{and}
 \,\,\,\,\,
 [\mathcal{L}_{\bf a}: A] = p,
    \end{equation}
    where $\mathcal{L}_{\bf a}$ is defined in (\ref{ppil}).

    Then $A = \Lambda$.}

     \vskip+0.3cm
     Proof.  The proof follows directly step-by-step to the proof of Lemma 1 and its Corollary from \cite{MOS}.
     The only minor changes are that
     one should consider integer parameters $M_1$ and $M_2$ instead of $2$ and $ n+1$ and that form
     $\varphi $ may be not integer. 
     $\Box$
     
     From Lemma 3  and Lemma 2 from \cite{MOS} we immediately deduce 
     \vskip+0.3cm
     {\bf Corollary.}
     {\it The maximal number of pairwise equivalent forms of the form 
     $Q_{\bf a}$ with ${\bf a}\in \mathscr{W}$ satisfying (\ref{ppiV}) is at most the number of integer  automorphisms of  the form $Q$.}

        \vskip+0.3cm
        
       {\bf 7. Proof of Theorem 2.}
         \vskip+0.3cm
       
       The proof follows the argument from Sections 5,6 from   \cite{MOS}. 
       For $n$ large enough take prime $p $ in the range
       \begin{equation}\label{priii}
       \frac{H_n}{2} \le p \le H_n,\,\,\,\,\,
       \text{where}
       \,\,\,\,\,
       H_n = \frac{27}{64} \, \frac{\Gamma\!\left(\tfrac{n}{2}+1\right)\sqrt{n+7}}{2^{{3n}/{2}}\pi^{n/2}} =
       \exp\left(\frac{n\log n}{2} - O(n)\right).
       \end{equation}
       We may assume that $ p>3(n+7)$  and so 
       \begin{equation}\label{pipip}
       (p,6(n+7)) = 1.
       \end{equation}
       Consider the set 
       $$ 
       \mathscr{W}' = \{ {\bf a}\in \mathscr{W}:\,\,\,
       \mathcal{L}_{\bf a} \cap \frak{A}_{f_n} \mathcal{E}_{f_n} =\Lambda_{f_n} \cap  \frak{A}_{f_n} \mathcal{E}_{f_n}
         =     S_{f_n}' \},
         \,\,\,\,\,
         \text{where}
         \,\,\,\,\,
         S_{f_n}'=
         S_{f_n}^0\cup \{{\bf 0}\},\,\,\,\,\,
          \frak{A}_{f_n} \mathcal{E}_{f_n} = \{{\bf z}\in \mathbb{R}^n:\,\, |{\bf z|}^2 \le 2\},
 $$
 and $S_{f_n}^{0}$ is as in (\ref{seet}).
For the number of elements of this set we have lower bound 
         \begin{equation}\label{lb}
       |\mathscr{W'}|
        \ge \frac{|\mathscr{W}|}{2} = \frac{p^{n-1}}{2} \ge \exp\left(\frac{n^2\log n}{2} - O(n^2)\right)
        \end{equation}
         because of the lower bound for $p$ from (\ref{priii}) and  the chain of inequalities  
       $$
      2\,     (|\mathscr{W}| -|\mathscr{W}'|) \le
       \sum_{{\bf a}\in \mathscr{W}} \left| \mathcal{L}_{\bf a} \cap (  \frak{A}_{f_n} \mathcal{E}_{f_n}\setminus S_{f_n}')\right|
       < \left|
       \frac{1}{p} \Lambda_{f_n}\cap   \frak{A}_{f_n} \mathcal{E}_{f_n}\right|  
      \le
      { (p+2n)^n \omega_n}=
      \frac{8}{3\sqrt{3}}\, 
      \frac{ (p+2n)^n (2\pi)^{n/2}}{\sqrt{n+7}\,\Gamma\!\left(\tfrac{n}{2}+1\right)} \le p^{n-1} = |\mathscr{W}|
       $$
       (we used equalities from (\ref{wo}) and right inequality from (\ref{priii})).
 
  \vskip+0.3cm
  
  Theorem 2 follows from the observation that 
  
    \vskip+0.3cm
    \noindent
  1) for  every $ {\bf a}\in \mathscr{W}'$  
  for quadratic form  
  \begin{equation}\label{fnba}
   f_{n,{\bf a}} = f_n(\frak{B}_{\bf a} {\bf x})
   \end{equation}
  we have
  $S_{   f_{n,{\bf a}} } = S_{   f_n }$,  and so  by Lemma 2 every $   f_{n,{\bf a}}$ is perfect but non-eutactic;
 
     \vskip+0.3cm
    \noindent
  2)    
  conditions (\ref{dudud}) and (\ref{pipip}) allow to apply Lemma  3 and its 
   Corollary together with   inequality (\ref{aut}), and in such a way we show that the number of pairwise equivalent forms of the form 
  (\ref{fnba})  with condition (\ref{ppiV})  with $ \varphi({\bf x}) = f_n ({\bf x})$ is $ \exp\left(O( n\log n)\right)$;

       \vskip+0.3cm
    \noindent
  3)   we have lower bound (\ref{lb}),
    meanwhile the number of vectors ${\bf a}\in \mathscr{W}$
    not satisfying  (\ref{ppiV})
      with $ \varphi({\bf x}) = f_n ({\bf x})$  
  is $ O(p^{n-2})$.
  
       \vskip+0.3cm
       Theorem 2 is proven.$\Box$

\vskip+0.3cm
{\bf Acknowledgements.} This paper is motivated by a question by Prof. Lenny Fukshansky about the number of eutactic forms in $n$ variables during the author's lecture at the conference 
"Diophantine Approximation and Related Topics",
Urgench State University, Urgench, Uzbekistan
August 24 - 28, 2026,
https://sites.google.com/view/urgench-diophantine-2026/home.
The author is very grateful to him.

The research is supported by   Austrian Science Fund (FWF), Forschungsprojekt PAT1961524.

\vskip+0.3cm
{\bf AI declaration.}  In preparation of this article the author used Anthropic's Claude   tools.

AI was used to verify all the proofs, correct issues and typos and  to search for the literature. An important idea suggested by  AI was  to use triangulation  of the fundamental domain of Voronoï polyhedron to adapt van Woerden's argument to prove inequality (\ref{vory2}). As for Theorem 2, the proofs were produced 'by hands' just generalising the ideas from \cite{MOS}, AI was used here to make calculations with Barnes' forms and to verify all the exposition.

 \end{document}